\documentclass[11pt,reqno]{amsart}

\usepackage[utf8]{inputenc}

\usepackage{amsmath,amsthm,amssymb}
\usepackage{amssymb}

\usepackage{graphics}
\usepackage{hyperref}
\usepackage[usenames, dvipsnames]{xcolor}

\usepackage{mathtools}
\mathtoolsset{showonlyrefs}

\usepackage[square,sort,comma,numbers]{natbib}

\definecolor{darkblue}{rgb}{0.0,0.0,0.3}
\hypersetup{colorlinks,breaklinks,
  linkcolor=darkblue,urlcolor=darkblue,
anchorcolor=darkblue,citecolor=darkblue}

\usepackage{booktabs}
\usepackage{threeparttable}

\theoremstyle{plain}
\newtheorem{theorem}{Theorem}
\newtheorem*{theorem*}{Theorem}

\newtheorem*{proposition*}{Proposition}

\newtheorem*{corollary*}{Corollary}

\newtheorem*{conjecture*}{Conjecture}

\theoremstyle{definition}
\newtheorem{remark}[theorem]{Remark}

\newtheorem{example}[theorem]{Example}

\theoremstyle{plain}
\newtheorem{prop}{Proposition}

\theoremstyle{definition}

\theoremstyle{remark}

\newcommand{\bP}{\mathbb{P}}
\newcommand{\bQ}{\mathbb{Q}}

\newcommand{\wB}{\widetilde{B}}
\newcommand{\wS}{\widetilde{S}}
\newcommand{\wX}{\widetilde{X}}
\newcommand{\wY}{\widetilde{Y}}

\newcommand{\Bl}{\operatorname{Bl}}

\newcommand{\ra}{\rightarrow}

\makeatletter
\let\@wraptoccontribs\wraptoccontribs
\makeatother

\title{Congruent Number Triangles with the Same Hypotenuse}
\author[David Lowry-Duda]{David Lowry-Duda}
\contrib[With an appendix by]{Brendan Hassett}
\thanks{The author was supported by the Simons Collaboration in Arithmetic
Geometry, Number Theory, and Computation via the Simons Foundation grant 546235.}
\thanks{The author wants to thank John Cremona, Brendan Hassett,
Joe Silverman, and Damiano Testa for the encouraging remarks and discussion.
The author also wants to thank an anonymous Referee for several comments, and
for directing the author to the paper~\cite{watkins2014ranks}.}
\thanks{BH was supported by National Science Foundation grant 1701659
and Simons Foundation award 546235}

\begin{document}

\maketitle

\begin{abstract}
In this article, we discuss whether a single congruent number $t$ can have two
(or more) distinct corresponding triangles with the same hypotenuse. We
describe and carry out computational experimentation providing evidence that
this does not occur.
\end{abstract}

\section{Introduction and Motivation}

A \emph{congruent} number is a rational number that appears as the area of a
rational right triangle. The \emph{congruent number classification problem} is
the problem of determining whether a given number $t$ is congruent. By scaling
the triangle, one may consider only the classification of squarefree integers
$t$; we do this throughout the article. The study of congruent numbers is
ancient and has a long history; see the survey article of Conrad~\cite{conrad}
for more.

We consider whether a single congruent number $t$ can have two
distinct rational right triangles with the same hypotenuse. This question was
first motivated by an observation in~\cite{hkldw_congruent}.

In that paper, Hulse, Kuan, Walker, and the author consider a convolution sum
whose main term exists only when $t$ is congruent. Let $\tau(n)$ denote the
square indicator function, taking $1$ if $n$ is a square and $0$ otherwise.
The primary theorem of~\cite{hkldw_congruent} is that
\begin{equation}
  \sum_{n \leq X} \sum_{m \leq X}
  \tau(m+n) \tau(m) \tau(m-n) \tau(nt)
  =
  C_t \sqrt X + O_t({(\log X)}^{r/2}),
\end{equation}
where $r$ is the rank of the elliptic curve $E_t: Y^2 = X^3 - t^2 X$ and
\begin{equation}
  C_t = \sum_{h \in \mathcal{H}(t)} \frac{1}{h}
\end{equation}
is the (convergent) sum over the set of hypotenuses $\mathcal{H}(t)$ of
dissimilar primitive right integer triangles with area equal to $t u^2$ (i.e.\
the squarefree part of the area is $t$).

Although there is implicit dependence on $t$ in the error term and poorly
understood dependence on $t$ in the main term, heuristically one should expect
that for large $X$, the main term ``quickly'' dominates.

In~\cite{hkldw_congruent}, it is also shown that the set of hypotenuses
$\mathcal{H}(t)$ is only logarithmically dense in $\mathbb{Z}_{\geq 1}$.
Thus heuristically the coefficient $C_t$ of the main term is well-estimated
by its first term $m/h$, where $h$ is the smallest hypotenuse in
$\mathcal{H}(t)$ and $m$ is the multiplicity of this hypotenuse. Further, the
main term should heuristically grow larger than the error term quickly once $X
> h/m$.

The question of the size of the smallest hypotenuse in $\mathcal{H}(t)$ is
closely related to the least height of a non-torsion point on $E_t$ (we make this
more precise in~\S\ref{sec:method}), but the multiplicity $m$ is more
mysterious.

In this article, we describe numerical experimentation suggesting that the
multiplicity is always exactly $1$. We also propose that this is always the case
as a conjecture.

\begin{conjecture*}
  There do not exist two dissimilar primitive right triangles with the same
  hypotenuse and whose squarefree parts of the areas are equal to each other.
\end{conjecture*}

In~\S\ref{sec:algebra}, we describe why this problem might be hard to fully
resolve. Hassett describes this more completely in the Appendix.

We then describe and carry out numerical experimentation using the free and open
source math software SageMath~\cite{sage}.
Our experimentation comes in two forms. Firstly, we do a deep
investigation for those congruent numbers below $1000$. For each congruent
number, we generate many different triangles by generating rational points
$E_t(\mathbb{Q})$ and the corresponding right triangles.

Secondly, we do a broad investigation for several congruent numbers
corresponding to curves $E_t$ of high rank. These curves will have more
rational points up to a bounding height, heuristically leading to more
triangles with hypotenuses up to some bound. In this investigation, we make
heavy use of~\cite{watkins2014ranks}, which gives several congruent number
curves $E_t$ of rank $6$ and $7$. In total, we investigate $1513$ curves of
ranks $6$ and $7$.

We describe these approaches and the results in~\S\ref{sec:method} and
in~\S\ref{sec:results}, respectively.

\section{Algebraic Formulation}\label{sec:algebra}

One can parametrize right triangles with two variables
$(s, t)$ such that the sides are given by $(s^2 - t^2, 2st, s^2 +
t^2)$. We follow the convention that the hypotenuse is written last in any
triple $(a, b, c)$ giving a right triangle, and that $s > t$. In this
correspondence, primitive integer right triangles correspond to relatively prime
integers $s$ and $t$.

Thus finding two right triangles with the same hypotenuse and whose areas are
the same up to multiplication by squares can be reformulated as finding
integers $(s, t, S, T)$ such that
\begin{equation}
  s^2 + t^2 = S^2 + T^2
\end{equation}
and where
\begin{equation}
  \frac{st(s^2 - t^2)}{ST(S^2 - T^2)} \; \text{is a square}.
\end{equation}
This last equation can be rewritten as
\begin{equation}
  u^2 st(s^2 - t^2) = v^2 ST(S^2 - T^2)
\end{equation}
for some positive integers $u, v$.

These two equations define a surface
\begin{equation}
  X \subset \mathbb{P}^3_{[s, t, S, T]} \times \mathbb{P}^1_{[u, v]}.
\end{equation}
If we could understand all points on this surface, we could almost certainly
resolve the conjecture. Unfortunately, this understanding appears to be beyond
current algebraic techniques. This is considered by Hassett in the Appendix,
where it is shown that the surface $X$ admits a resolution of singularities
that is of general type and simply connected (see
Proposition~\ref{prop:appendix}).

Thus we expect that a complete resolution of the conjecture is not
within reach.

\section{Description of Methodology}\label{sec:method}

With a purely theoretical resolution likely out of reach, we turn now to
computational experimentation. There is a well-known correspondence between
right triangles $(a, b, c)$ with $ab/2 = t$ (which we may assume without loss of
generality is a squarefree integer) and $\mathbb{Q}$-rational points on the
elliptic curve $E_t: Y^2 = X^3 - t^2 X$ where $Y \neq 0$.
The inverse maps of this correspondence are given by
\begin{align}
  (a, b, c)
  &\mapsto \Big( \frac{tb}{c-a}, \frac{2t^2}{c-a} \Big)
  = (X, Y),
  \\
  (X, Y)
  &\mapsto \Big( \frac{X^2 - t^2}{Y}, \frac{2tX}{Y}, \frac{X^2 + t^2}{Y} \Big)
  = (a, b, c).
\end{align}
(See~\cite{koblitz} for a historical overview and description of the relationship
between congruent numbers and elliptic curves).

\subsection{Enumeration through elliptic curves}

This provides an explicit and computable way of enumerating all rational right
triangles with area $t$: find generators for $E_t(\mathbb{Q})$ and enumerate
their linear combinations. This method of enumeration is of course vastly
superior to naive enumeration and is at the core of the organization of the
numerical experimentation.

In principle, it is necessary to also include the torsion subgroup of
$E_t(\mathbb{Q})$. However for these curves, the torsion subgroup is completely
understood and can be ignored. In particular, the torsion subgroup is isomorphic
to $\mathbb{Z}/2\mathbb{Z} \times \mathbb{Z}/2\mathbb{Z}$, generated by two
points $T_1$ and $T_2$. If the triangle $(a, b, c)$ corresponds to the point $P
\in E_t(\mathbb{Q})$, the eight equivalent triangles $(\pm a, \pm b, \pm c)$
correspond to the points $\pm P + \epsilon_1 T_1 + \epsilon_2 T_2$, where
$\epsilon_i \in \{0, 1\}$. Thus we omit no triangles (up to similarity) by
omitting consideration of torsion points.

\subsection{Hypotenuses, and heights, grow exponentially}%
\label{ssec:hypotenuses_and_heights}

Given an elliptic curve with known rank and generators, the primary obstruction
to collecting numerical evidence is the size of the triangles and hypotenuses.
There is a general theory that the heights of points on elliptic curves grow
very rapidly under repeated addition. In Corollary~1 of~\cite{hkldw_congruent},
it is shown that
\begin{equation}
  \lvert \{ h \in \mathcal{H}(t) : h \leq  X/2t \} \rvert
  =
  O_t({(\log X)}^{r/2}),
\end{equation}
where $r$ is the rank of the corresponding elliptic curve $E_t$. In particular,
letting $H(P)$ denote the naive height of the $X$-coordinate of $P$, it is
shown that the hypotenuse $h$ corresponding to a point $P$ on $E_t$ satisfies $h
\geq H(P)/2t$, and thus hypotenuses grow at least as quickly as the heights.

\begin{example}
The elliptic curve $E_6: Y^2 = X^3 - 36 X$ has rank $1$ and the free part of
$E_6(\mathbb{Q})$ is generated by $g = (-3, 9)$. This point corresponds to the
primitive right triangle $(3, 4, 5)$. The hypotenuse is 3 bits long.
The point $2*g = (25/4, -35/8)$ corresponds to the primitive right triangle
$(49, 1200, 1201)$, whose hypotenuse is 11 bits long. The point $4*g$ has
hypotenuse $2094350404801$, which is 29 bits long. The point $400*g$ (the
largest multiple we considered on this curve) corresponds to a right triangle
whose hypotenuse requires 410426 bits to store.
\end{example}

As heights grow exponentially, we see that hypotenuses grow exponentially.
If we assume that the hypotenuses behave like exponentially growing random
variables, then heuristically we would expect that it is very unlikely that
two points of large height correspond to the same hypotenuse.
Informally, we should expect that if a pair of points yields two triangles
with hypotenuse that are close in size, then both points are probably of small
height.

For curves of rank $1$, it is impractical to compute more than a
few hundred triangles --- they simply grow too quickly. For curves of higher
rank, however, it is possible to compute millions of triangles before the
individual sides of the triangles become prohibitively large.

\subsection{Two strategies of experimentation}

We investigate two potential sources for a counterexample to the conjecture:
when all numbers are small and coincidental collisions are most likely, and
when there are many triangles up to a certain size due to high rank. Thus
we organize our experimentation into two strategies: one for elliptic curves
$E_t$ with small $t$, and one for high rank elliptic curves.

First, we study the curves $E_t$ for all squarefree $t \leq 1000$. Each of
these curves has rank $0$, $1$, or $2$. In total, this includes $361$ curves
corresponding to congruent numbers. These are the curves with the smallest $t$
and with the simplest triangles.

As these curves have low rank, we can carry out ``deep'' investigation on each
curve. More specifically, if $\{g_i\}$ denote generators for the free part of
$E_t(\mathbb{Q})$, then we can study triangles coming from $\sum a_i g_i$ for
all coefficients $\{a_i\}$ in a large box $\lvert a_i \rvert \leq K$ for a
constant $K = K(t)$. For rank $1$ curves, we choose $K(t) \geq 300$, and for
rank $2$ curves we choose $K(t) \geq 75$. As the points $P$ and $-P$ generate
the same triangle, this translates to considering at least $300$ triangles from
each rank $1$ curve and at least $(151^2 - 1)/ 2 = 11400$ triangles from each
rank $2$ curve.

Second, we study curves $E_t$ of known high rank. To get these curves, we use
the large-scale project detailed in~\cite{watkins2014ranks}, which builds off
of earlier work of Rogers~\cite{rogers2000rank,rogersthesis}. In this project,
Watkins et al.\ investigated ranks of congruent number curves. They identify
$1486$ curves $E_t$ of rank at least $6$ and $27$ curves of rank $7$, giving
$1513$ curves in total of high rank. We examine rational points on each of
these $1513$ curves.

We note that no curves of higher rank are known, and it is conjectured that
ranks in the family of congruent number curves are
bounded~\cite{honda1961isogenies}. In \S11 of \cite{watkins2014ranks}, it is
further demonstrated that a heuristic of Granville's suggests that the maximum
rank of a congruent number elliptic curve is $7$.

As noted in \S\ref{ssec:hypotenuses_and_heights}, we should heuristically
expect that it is more likely for points of small height to give two triangles
with the same hypotenuse. For this reason, we carry out a ``wide but
shallow'' investigation of these curves. Again letting $\{g_i\}$ denote a set of
generators for the free part of $E_t(\mathbb{Q})$, we study triangles coming
from points $\sum a_i g_i$, where all coefficients $\{ a_i \}$ are in the small
box $\lvert a_i \rvert \leq 4$. This gives $267520$ triangles from each rank
$6$ curve and $2391484$ triangles from each rank $7$ curve.

\begin{remark}
  It is likely that using the box $\lvert a_i \rvert \leq 2$ would be
  sufficient. In \S\ref{sec:results}, we observe that $\lvert a_i \rvert \leq
  2$ finds all closest pairs of hypotenuses on the curves investigated. We use
  the larger box to conduct a more robust investigation.
  Computing and manipulating these points is by far the most time intensive
  portion of this experiment. With the larger coefficient boxes $\lvert a_i
  \rvert \leq 4$, we spend approximately $32$ minutes on average for each rank
  $6$ curve, and approximately $5$ hours on average for each rank $7$ curve. In
  total, we used approximately $940$ CPU hours on this part of the
  experimentation.
%
\end{remark}

\subsection{Computing rank and generators}

In order to generate triangles for a congruent number $t$, we require that we
can compute generators for (the free part of) $E_t(\mathbb{Q})$. We perform
computations for two different sets of curves: those curves $E_t$ with $t <
1000$ and for several chosen curves of higher rank.

For the $1513$ curves of higher rank, we use ranks and generators available
as an electronic supplement to~\cite{watkins2014ranks}.\footnote{The supplement
to their paper made this data available in Magma code. We've converted it to a
SageMath-friendly format and made it available at
\url{https://github.com/davidlowryduda/notebooks/blob/master/Papers/largetdata.sage}}
We note that not all of these curves have explicitly known rank: all have rank
at least $6$, and $27$ have rank $7$, but Watkins et al.\ couldn't verify the
ranks of each curve. We do not try to complete this verification here, and
instead use the $6$ known generators for each curve of presumed rank $6$.

Curves corresponding to $t < 1000$ can be handled within SageMath~\cite{sage}.
Directly using SageMath's \texttt{EllipticCurve.rank} and
\texttt{EllipticCurve.gens} works for all but $84$ cases. The smallest
problematic example is $113$. We note that the curve $E_{113}$
is contained within the L-function and Modular Form Database
(LMFDB)~\cite{lmfdb}. However, it is also possible to directly use Tunnell's
criterion~\cite{tunnell} to quickly confirm that $113$ is not congruent.
Applying Tunnell's criterion on the remaining $83$ cases to remove
non-congruent numbers reduces this to $54$ cases. The smallest remaining
example is $157$.

The number $157$ is congruent and is somewhat famous for its late
classification by Zagier~\cite{zagier1989elliptische}; the simplest rational
triangle with area $157$ is extremely complicated.

The techniques SageMath uses to compute the rank and generators for elliptic
curves are based on John Cremona's MWRANK~\cite{john_cremona_2019_2577868}).
SageMath includes a copy of MWRANK and a partial interface to MWRANK's
functionality. Calling MWRANK directly and allowing up to 2000 seconds of
computation time per curve provides generators and ranks for an additional $20$
cases (including $157$). This leaves $34$ congruent numbers up to $1000$, the
smallest of which is $277$.

To understand the remaining $34$ curves, we use functionality from
PARI-GP~\cite{PARI}, which is also packaged with SageMath. Although MWRANK
struggles to exactly compute the rank and generators, it can compute that the
upper bound on the rank is $1$ for each of these curves. For the curve
$E_{277}$ for example, calling MWRANK through Sagemath with the command
\texttt{EllipticCurve([-277*277,0]).rank\_bound()} shows an upper bound of $1$
for the rank.

Thus one can use the Heegner point method~\cite{elkies1994heegner} to find a
generator for the rational points on these curves. The existence of such a
generator confirms that the rank is exactly $1$. This is implemented remarkably
well in PARI. Calling
\texttt{pari(EllipticCurve([-277*277,0])).ellheegner()} from SageMath quickly
returns a generator for $E_{277}(\mathbb{Q})$.


\begin{remark}
  The CAS Magma~\cite{MAGMA} is also capable of generating ranks and generators
  for these elliptic curves. An anonymous Referee noted that
  \texttt{Generators(EllipticCurve([-277*277,0]))} returns a set of generators
  almost immediately, and functions almost as quickly for all congruent numbers
  up to $1000$.
\end{remark}

\section{Experimental Results}\label{sec:results}

We performed the strategies outlined in~\S\ref{sec:method} for each squarefree
congruent number $t$ with $t \leq 1000$ and for the $1513$ curves of higher
rank from~\cite{watkins2014ranks}.
(See the author's demonstration github~\cite{dld_notebook_cnp_hypotenuse} for a
reference implementation to obtain and manipulate generators when $t \leq
1000$; for higher $t$, we do similar manipulations with the generators found by
Watkins et al).

\subsection{Results for $t \leq 1000$}

Out of the $608$ squarefree numbers up to $1000$, we verified that $327$ of
them are congruent numbers and produced several triangles for each one.
Of these, $274$ corresponding to curves with rank 1 and $53$ to curves with
rank $2$. We note that this is consistent with the celebrated conjecture
of Goldfeld~\cite{goldfeld1979conjectures}, which implies that although there
may be infinitely many such curves of rank at least $2$, these should be sparse
and correspond to $0$\% of all congruent numbers in natural density.

Among all these computed triangles, we found no counterexample to the
conjecture. There were no two dissimilar right triangles with the same
hypotenuse (after scaling to a primitive right triangle) corresponding to the
same congruent number.

Recall that we computed at least $300$ triangles from each curve of rank $1$
and at least $11400$ from each curve of rank $2$. Per the heuristic in
\S\ref{ssec:hypotenuses_and_heights}, we would expect that pairs of points
corresponding to triangles with the close hypotenuses should come from points
of small height.

This is true for all the curves tested for $t \leq 1000$. Further, the nearest
pairs of hypotnuses among triangles associated to the same curve were always
among the smallest hypotenuses. In terms of our generators $\{ g_i \}$, the
smallest pair of hypotenuses always came from triangles corresponding to the
points $ \sum a_i g_i $ with $\lvert a_i \rvert \leq 2$. (Although there are
not cannonical choices of generators, we note that MWRANK returns generators
with small naive height~\cite{cremona1997algorithms}).

\begin{example}
On $E_6(\mathbb{Q})$, whose free part is generated by $(-3, 9)$, the two
nearest hypotenuses come from the triangles $(3,4,5)$ and $(49, 1200, 1201)$,
which come from $g$ and $2g$ respectively.

For $t = 34$, we have that $\mathrm{rank}\big(E_{34}(\mathbb{Q})\big) = 2$. The free
part of $E_{34}(\mathbb{Q})$ is generated by the points $g_1 = (-16, 120)$ and
$g_2 = (-2, 48)$. The two closest hypotenuses come from the triangles $(225,
272, 353)$ and $(17, 144, 145)$, corresponding to $g_1$ and $g_2$ respectively.
\end{example}

This adds support for our choice of box $\lvert a_i \rvert \leq 4$ for points
on higher rank elliptic curves.

Although no pair of hypotenuses exactly matched, we can ask about how close two
hypotenuses could be. Of those congruent numbers $t \leq 1000$, we determined
that $5$ have pairs of dissimilar triangles with hypotenuses that differ by
less than $t$: $210, 330, 546, 609, 915$. All five of these come from rank $2$
curves.

The curve $E_{210}$ has a particularly close pair of triangles whose
hypotenuses differ by $8$. This is the pair of triangles $(12, 35, 37)$ and
$(20, 21, 29)$, which correspond to the two generators for the free part of
$E_{210}(\mathbb{Q})$.

\subsection*{Smallest hypotenuses}

\begin{figure}[ht]
  \caption{%
  Smallest hypotenuses for congruent numbers $t \leq 1000$\label{fig:plot}
  }
  \includegraphics[scale=0.9]{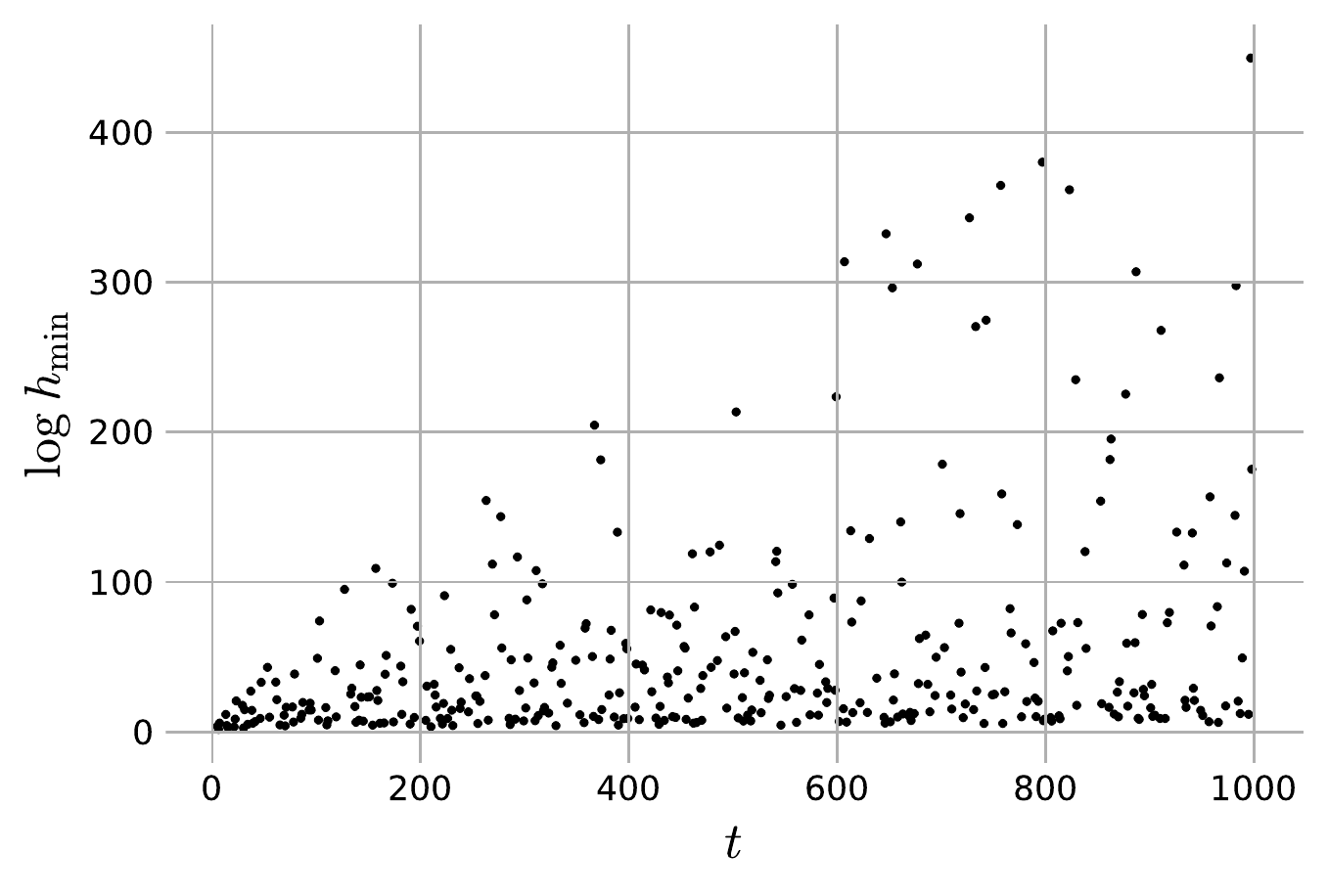}
  \centering
\end{figure}

Further, we took this opportunity to collect data on the size of the smallest
hypotenuse for each congruent number $\leq 1000$. In Figure~\ref{fig:plot}, we
plot squarefree congruent numbers $t$ against the log of the smallest
hypotenuse of a primitive right triangle with squarefree part of the area equal
to $t$.

Within the plot, the first congruent number $t$ with smallest hypotenuse greater
than $e^{100}$ is Zagier's number, $157$. The \emph{largest} smallest hypotenuse
we computed corresponds to the congruent number $t = 997$, and is approximately
$e^{449}$.

Many smallest-hypotenuses appear to be moderate in size, but no consistent
pattern emerges.

\subsection{Results for curves of rank $6$ and $7$}

We produced hundreds of thousands of triangles for each of the $1513$ curves of
high rank. We found no counterexample to the conjecture.

As above, we tracked closest pairs of hypotenuses for each curve. Even though
we computed all rational points $\sum a_i g_i$ in the coefficient box where
each coefficient satisfies the bound $\lvert a_i \rvert \leq 4$, we again
observed that for each $t$, the closest pair of hypotenuses came from points in
the smaller box $\lvert a_i \rvert \leq 2$. The heuristic from
\S\ref{ssec:hypotenuses_and_heights} appears to continue to hold, indicating
that nearest pairs of hypotenuses occur among the points of the least height on
the curve.

We consider how close hypotenuses from two different triangles could be.
Of all the $1513$ higher rank curves investigated, $375$ has two hypotenuses
that differed by less than the corresponding congruent number $t$.

In absolute terms, the closest pair came from $t = 6611719866 \approx 6.611
\times 10^9$, which has triangles with hypotenuses $30544225$ and $67119265$,
differing by about $3.657 \times 10^7$. This is very large.
We also noted a pair for the congruent number $t = 1902736244939034 \approx
1.9027 \times 10^{15}$, which has triangles with hypotenuses given approximately
by $1.061 \times 10^{10}$ and $1.109 \times 10^{10}$, differing by about $4.769
\times 10^8$. These are also the two smallest hypotenuses of triangles for $t$.
Although this is also very large, this is significantly smaller than both $t$
and the smallest hypotenuse.

Both of these ``near-miss'' examples came from rank $6$ elliptic curves.

\begin{figure}[ht]
  \caption{%
  Smallest hypotenuses from (presumed) rank $6$ curves\label{fig:plot2}
  }
  \includegraphics[scale=0.9]{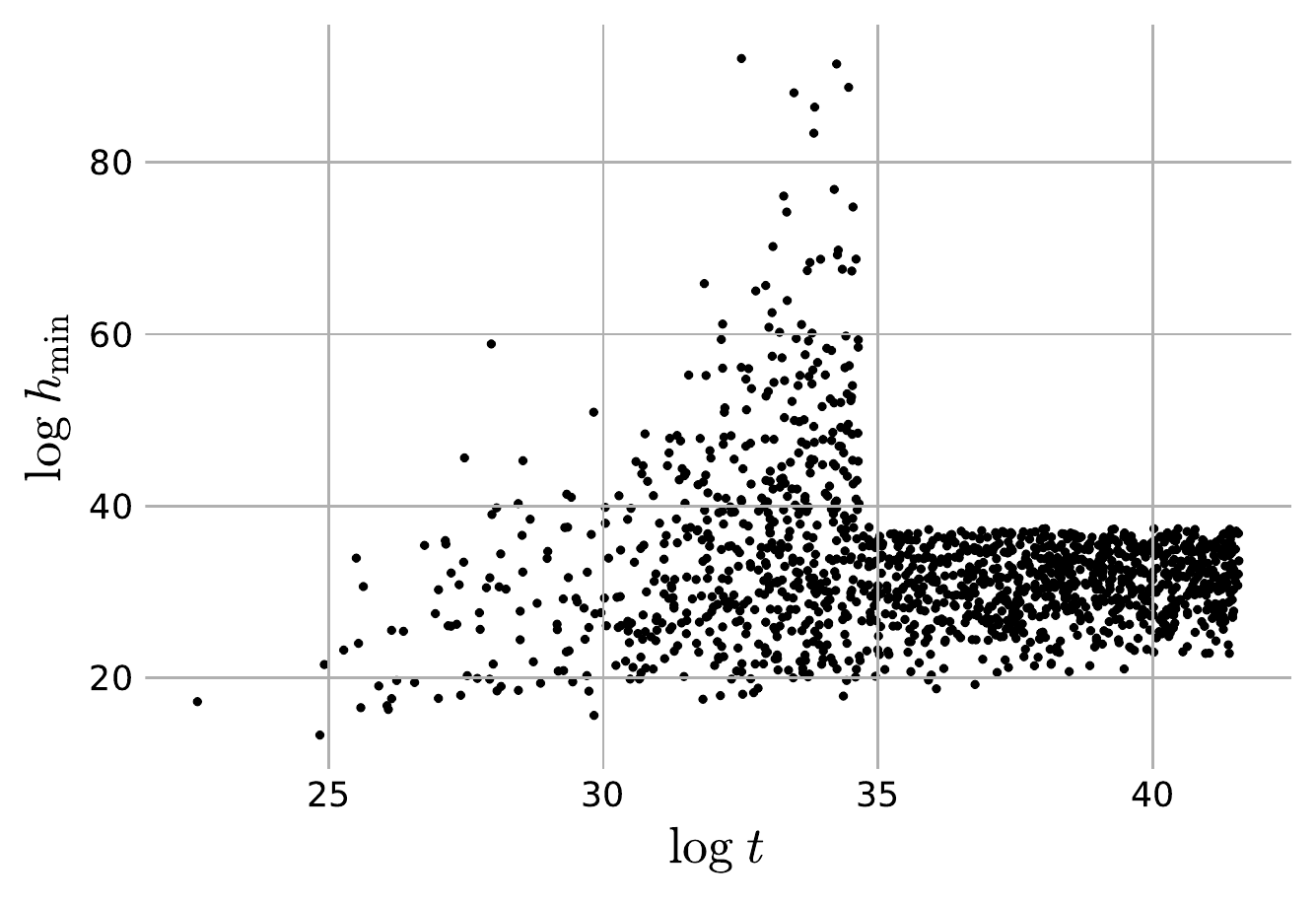}
  \centering
\end{figure}

\subsection*{Smallest hypotenuses}

We again took this opportunity to collect data on the size of the smallest
hypotenuse coming from each of the $1513$ curves of higher rank. As the numbers
are encompass a vastly larger domain, we present this in two $\log$-$\log$ plots.

In Figure~\ref{fig:plot2}, we show smallest hypotenuses for each $t$ from
curves of (presumed) rank $6$. We first explain the lack of points in the
upper-right portion of the plot. This is caused by the methodology used to
generate the $1513$ curves in~\cite{watkins2014ranks} --- they performed an
almost-exhaustive search for high rank elliptic curves for all $t \leq 2^{50}$,
while for higher $t$ they only looked for curves with points of relatively low
height. As $2^{50} \approx e^{34.5}$ and as the hypotenuse $h$ corresponding to
a point $P$ satisfies the bound $h > H(p)/2t$, this explains the missing
region of data. (Although we found it surprising that this plot reveals the
choices of bounds used to find the elliptic curves so clearly).

\begin{figure}[htb]
  \caption{%
  Smallest hypotenuses from rank $7$ curves\label{fig:plot3}
  }
  \includegraphics[scale=0.9]{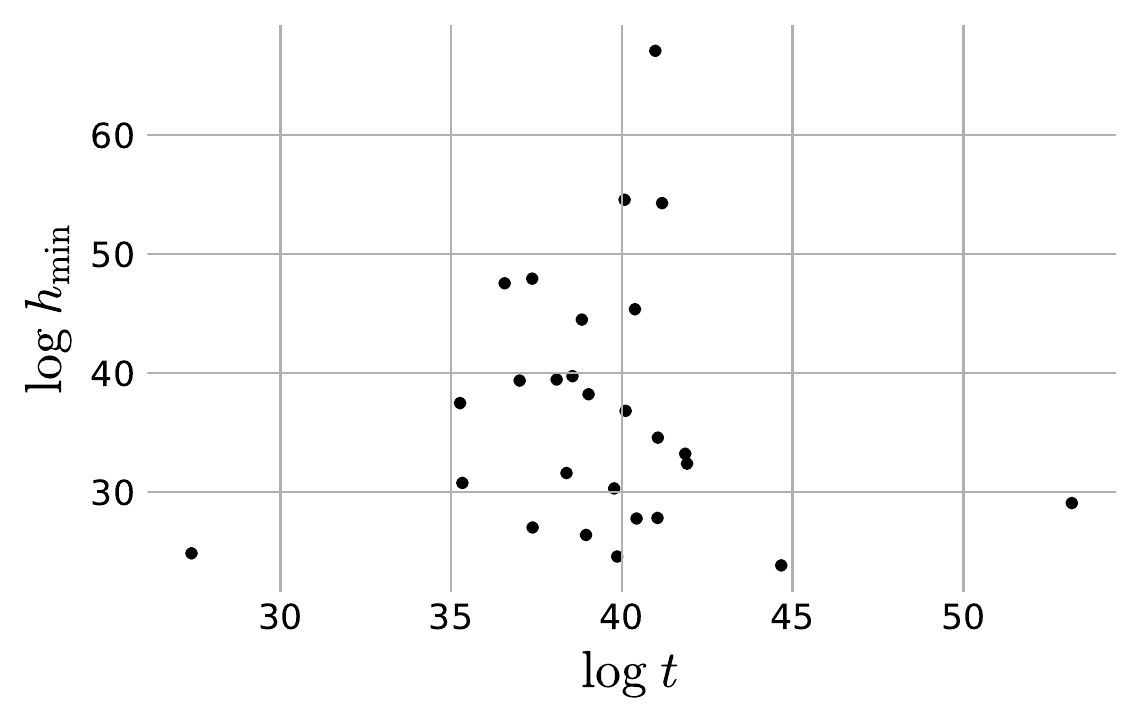}
  \centering
\end{figure}

Comparing Figure~\ref{fig:plot2} with Figure~\ref{fig:plot}, we note the
disparity in size. Every hypotenuse coming from a higher rank elliptic curve
was smaller than than the smallest hypotenuse on $E_{157}$.
On the other hand, each smallest hypotenuse was larger than $e^{15}$.

In Figure~\ref{fig:plot3}, we show smallest hypotenuses coming from curves of
rank $7$. As with the previous plot, there are patterns among the curves that
result from the way in which the curves were initially found.
In~\cite{watkins2014ranks}, they performed an extensive search for curves of
rank $7$ for $t \leq 2^{60} \approx e^{41.6}$, and for higher $t$ they only
looked for curves with rational points of low height. Although it is less
obvious in this plot, we note that we should expect that there is a lot of
missing data in the upper right portion of this plot.

%
%

\appendix

\section{Invariants of the Lowry-Duda moduli space}
\begin{center}by \textsc{Brendan Hassett}\end{center}

\renewcommand*{\theprop}{\Alph{prop}}

We verify that the moduli space $X$ introduced by Lowry-Duda admits a resolution of singularities
$\wX \ra X$ that is of general type and simply connected. Since $\wX(\bQ) \neq \emptyset$, existing
Diophantine techniques (like Faltings' Theorem \cite{Fal}) shed little light on
the structure of the rational points on $X$.

\subsection*{Double cover realization} Let
$$X\subset \bP^3_{[s,t,S,T]} \times \bP^1_{[u,v]}$$
be the surface over $\bQ$ defined by the equations
\begin{align*}
s^2+t^2&=S^2+T^2 \\
u^2st(s^2-t^2)&=v^2ST(S^2-T^2).
\end{align*}
The first equation defines a surface $Y\subset \bP^3$
isomorphic to
$\bP^1 \times \bP^1$.
Projection onto the first factor induces a morphism
$\varpi: X \ra Y$
that is generically finite of degree two. The branch curve
$$B:=\{st(s^2-t^2)ST(S^2-T^2)=0 \} \subset Y$$
is a union of eight planar sections of $Y$, each a smooth conic curve defined over $\bQ$.

The morphism $\varpi$ fails to be flat
over the locus
$$
Z:=\{st(s^2-t^2)=ST(S^2-T^2)=0\} \subset B \subset Y
$$
consisting of the following 32 points:
\begin{align*}
[0,1,0,\pm 1] &= \{s=0\} \cap \{S=0\} \\
[0, 1, \pm 1, 0] &= \{s=0\} \cap \{T=0\} \\
[0,\pm \sqrt{2}, 1,1] &= \{s=0\} \cap \{S-T=0\} \\
[0,\pm \sqrt{2}, 1,-1] &= \{s=0\} \cap \{S+T=0\} \\
[1,0,0,\pm 1] &= \{t=0\} \cap \{S=0\} \\
[1,0,\pm 1,0] &= \{t=0\} \cap \{T=0\} \\
[\pm \sqrt{2},0, 1,1] &= \{t=0\} \cap \{S-T=0\} \\
[\pm \sqrt{2},0, 1,-1] &= \{t=0\} \cap \{S+T=0\} \\
[1,1,0,\pm \sqrt{2}] &= \{s-t=0\} \cap \{S=0\} \\
[1,1,\pm \sqrt{2},0] &= \{s-t=0\} \cap \{T=0\} \\
[1,1, \pm 1, \pm 1] &= \{s-t=0\} \cap \{S-T=0\} \\
[1,1, \pm 1,\mp 1] &= \{s-t=0\} \cap \{S+T=0\} \\
[1,-1,0,\pm \sqrt{2}] &= \{s+t=0\} \cap \{S=0\} \\
[1,-1,\pm \sqrt{2},0] &= \{s+t=0\} \cap \{T=0\} \\
[1,-1, \pm 1,\pm 1] &= \{s+t=0\} \cap \{S-T=0\} \\
[1,-1, \pm 1,\mp 1] &= \{s+t=0\} \cap \{S+T=0\}
\end{align*}
Points of $Z$ are nodes of $B$.
The double cover of $Y$ branched over $B$ has $A_1$ singularities
over these $32$ nodes which are automatically resolved in $X$.

The remaining points of intersection among components of $B$ are:
\begin{align*}
 [0,0,1,\pm i] &=
  \{s=0\}\cap \{t=0\} \cap \{s-t=0\} \cap \{s+t=0\}\\
 [1,\pm i,0,0] &=
 \{S=0\}\cap \{T=0\} \cap \{S-T=0\} \cap \{S+T=0\}
\end{align*}
Write
$$W:=\{[0,0,1,\pm i],[1,\pm i, 0, 0]\},$$
where $B$ has multiplicity-four singularities.
There are $\binom{8}{2}=28$ pairs of conics altogether so all
intersections are accounted for.

\subsection*{Resolving the bad singularities}
We briefly review background on surface singularities:
A double cover of a smooth surface branched along a curve has
ADE singularities if and only if the curve has
simple (ADE) singularities.
An isolated hypersurface singularity $s\in S$ is {\em Du Val} if
there exists a resolution $\beta: \wS \ra S$ such that
$K_{\wS}-\beta^*K_S$ is effective. Complex Du Val singularities
coincide with the ADE surface singularities \cite[\S 4.20]{Reid}.
Given a projective surface $S$ with Du Val singularities,
pluricanonical differentials on $S$ are the same as pluricanonical
differentials on a resolution $\wS$. In particular, if $K_S$
is ample then $\wS$ is of general type.

Returning to $X$, we seek to resolve the singularities that are
not Du Val, which are associated with the multiplicity four singularities of $B$
along $W$.
We observe that
$$W\subset \{s^2+t^2=S^2+T^2=0\},$$
i.e., the vertices of a cycle of four rational curves contained in $Y$, whose
complement is a torus. The blow-up
$$\wY:=\Bl_W(Y) \ra Y$$
is therefore toric with boundary consisting of an octagon of $\bP^1$'s
with self-intersections alternating between $-1$ and $-2$.
The induced blow-up
$$\beta: \wX \ra X$$
resolves the singularities of $X$.

\subsection*{Canonical class of $X$}
First, we choose a basis for the Picard group of $\wY$: Let
$f_1$ and $f_2$ be the rulings of $Y\simeq \bP^1 \times \bP^1$ and
$E_1,\ldots,E_4$ the exceptional divisors over $W$ so
$$\operatorname{Pic}(\wY)=\left<f_1,f_2,E_1,E_2,E_3,E_4\right>.$$
We order so that the boundary octagon has sides with classes
$$E_1,f_1-E_1-E_2,E_2,f_2-E_2-E_3,E_3,f_1-E_3-E_4,E_4, f_2-E_4-E_1.$$

The canonical class
$$K_{\wY}\equiv -2f_1-2f_2+E_1+E_2+E_3+E_4.$$
The proper transform $\wB$ of $B$ -- noting the multiplicity four singularities along $W$ --
has class
$$8(f_1+f_2)-4(E_1+E_2+E_3+E_4);$$
$\wB$ is disjoint from the four $(-2)$ curves in the octagon.
Blowing down these four curves induces a birational morphism
$$\wY \ra \Sigma$$
to a quartic del Pezzo surface with four $A_1$ singularities. This morphism is
induced by the linear series $\lvert 2f_1+2f_2 - E_1 - E_2 - E_3 - E_4\rvert$ corresponding to the
quadratic equations for $W$ modulo the defining equation for $Y$.

We refer the reader to \cite[IV.22]{BHPV} for a discussion of the invariants
of double covers of surfaces branched along curves with simple singularities.

The canonical class of $\wX$ is obtained by pulling back
$$K_{\wY}+\frac{1}{2}\wB \equiv 2(f_1+f_2)-(E_1+\ldots+E_4)$$
through the induced morphism $\wX \rightarrow \wY$.
The image of $\wY$ and $\wX$ under this linear series
equals $\Sigma$. We summarize this in the following proposition.

\begin{prop}\label{prop:appendix}
The surface $\wX$ is of general type with $K_{\wX}^2=8$, and the
linear series $|K_{\wX}|$ induces a morphism
$$\wX \ra \Sigma$$
that is generically finite of degree two.
\end{prop}
In particular, we find
$$\dim \Gamma(K_{\wX})= \dim \Gamma(\wY,K_{\wY}+\frac{1}{2}\wB)=5.$$

\subsection*{$\wX$ is simply connected}
To compute the other invariants of $\wX$, we use observations
following immediately from simultaneous resolution of
ADE singularities \cite{Brieskorn}: $\wX$ is deformation equivalent to a
double cover of $\wY$ branched along a {\em generic} member of
$$\lvert \wB \rvert= \lvert 8(f_1+f_2)-4(E_1+E_2+E_3+E_4) \rvert,$$
disjoint from the $(-2)$-curves.
Smoothing the $A_1$ singularities of $\Sigma$ and noting that the
branch divisor lies in $\lvert -4K_{\Sigma}\rvert$, we find that $\wX$
is also deformation equivalent to a double cover of a quartic del
Pezzo surface $\Sigma'$ branched over a divisor $A\in \lvert -4K_{\Sigma'}\rvert$.
Such surfaces are simply connected by the Lefschetz hyperplane
theorem, hence
$$\pi_1(\wX)=\{1\}.$$

\subsection*{The remaining invariants}
The branch curve $A$
satisfies
$$\deg K_A = 48$$
by adjunction on $\Sigma'$.
Thus we conclude
$$\chi(\wX) = 2(\chi(\Sigma')-\chi(A))+\chi(A)=2(8+48)-48=64,$$
$b_2(\wX)=62$, and $h^1(\Omega^1_{\wX})=52$.

\vspace{20 mm}
\bibliographystyle{alpha}
\bibliography{bibfile}

\newcommand{\etalchar}[1]{$^{#1}$}
\begin{thebibliography}{BHPVdV04}

\bibitem[BBB{\etalchar{+}}00]{PARI}
Christian Batut, Karim Belabas, Dominique Bernardi, Henri Cohen, and Michel
  Olivier.
\newblock {\em User's Guide to PARI-GP}.
\newblock Universit{\'e} de Bordeaux I, 2000.

\bibitem[BCP97]{MAGMA}
Wieb Bosma, John Cannon, and Catherine Playoust.
\newblock The {M}agma algebra system. {I}. {T}he user language.
\newblock {\em J. Symbolic Comput.}, 24(3-4):235--265, 1997.
\newblock Computational algebra and number theory (London, 1993).

\bibitem[BHPVdV04]{BHPV}
Wolf~P. Barth, Klaus Hulek, Chris A.~M. Peters, and Antonius Van~de Ven.
\newblock {\em Compact complex surfaces}, volume~4 of {\em Ergebnisse der
  Mathematik und ihrer Grenzgebiete. 3. Folge. A Series of Modern Surveys in
  Mathematics [Results in Mathematics and Related Areas. 3rd Series. A Series
  of Modern Surveys in Mathematics]}.
\newblock Springer-Verlag, Berlin, second edition, 2004.

\bibitem[Bri68]{Brieskorn}
Egbert Brieskorn.
\newblock Die {A}ufl\"{o}sung der rationalen {S}ingularit\"{a}ten holomorpher
  {A}bbildungen.
\newblock {\em Math. Ann.}, 178:255--270, 1968.

\bibitem[CMP{\etalchar{+}}19]{john_cremona_2019_2577868}
John Cremona, Marcus Mo, Julien Puydt, Fernando Perez, François Bissey, Isuru
  Fernando, William Stein, qed777, Giovanni Mascellani, Keno Fischer, Jerry
  James, Dima Pasechnik, abergeron, and Antonio Rojas.
\newblock {\em JohnCremona/eclib v20190226}, February 2019.
\newblock Also at {\tt http://homepages.warwick.ac.uk/~masgaj/mwrank/}.

\bibitem[Con08]{conrad}
Keith Conrad.
\newblock The congruent number problem.
\newblock {\em Harvard College Mathematical Review}, 2:58--74, 2008.
\newblock http://www.math.harvard.edu/hcmr/issues/2a.pdf.

\bibitem[Cre97]{cremona1997algorithms}
John~E Cremona.
\newblock {\em Algorithms for Modular Elliptic Curves}.
\newblock Cambridge University Press, 1997.

\bibitem[Elk94]{elkies1994heegner}
Noam~D Elkies.
\newblock Heegner point computations.
\newblock In {\em International Algorithmic Number Theory Symposium}, pages
  122--133. Springer, 1994.

\bibitem[Fal94]{Fal}
Gerd Faltings.
\newblock The general case of {S}. {L}ang's conjecture.
\newblock In {\em Barsotti {S}ymposium in {A}lgebraic {G}eometry ({A}bano
  {T}erme, 1991)}, volume~15 of {\em Perspect. Math.}, pages 175--182. Academic
  Press, San Diego, CA, 1994.

\bibitem[Gol79]{goldfeld1979conjectures}
Dorian Goldfeld.
\newblock Conjectures on elliptic curves over quadratic fields.
\newblock In {\em Number Theory Carbondale 1979}, pages 108--118. Springer,
  1979.

\bibitem[HKLDW18]{hkldw_congruent}
Thomas~A. Hulse, Chan~Ieong Kuan, David Lowry-Duda, and Alexander Walker.
\newblock A shifted sum for the congruent number problem, 2018.

\bibitem[Hon61]{honda1961isogenies}
Taira Honda.
\newblock Isogenies, rational points and section points of group varieties.
\newblock In {\em Japanese journal of mathematics: transactions and abstracts},
  volume~30, pages 84--101. The Mathematical Society of Japan, 1961.

\bibitem[Kob93]{koblitz}
Neal Koblitz.
\newblock {\em Introduction to elliptic curves and modular forms}, volume~97 of
  {\em Graduate Texts in Mathematics}.
\newblock Springer-Verlag, New York, second edition, 1993.
\newblock http://dx.doi.org/10.1007/978-1-4612-0909-6.

\bibitem[LD]{dld_notebook_cnp_hypotenuse}
David Lowry-Duda.
\newblock Sage implementation for congruent number triangles with the same
  hypotenuse.
\newblock
  \url{https://github.com/davidlowryduda/notebooks/blob/master/Papers/CongruentNumberTrianglesSameHypotenuse.sage}.
\newblock Accessed: 2020-01-31.

\bibitem[{LMF}19]{lmfdb}
The {LMFDB Collaboration}.
\newblock The {L}-functions and modular forms database.
\newblock \url{http://www.lmfdb.org}, 2019.
\newblock [Online; accessed 30 October 2019].

\bibitem[Rei97]{Reid}
Miles Reid.
\newblock Chapters on algebraic surfaces.
\newblock In {\em Complex algebraic geometry ({P}ark {C}ity, {UT}, 1993)},
  volume~3 of {\em IAS/Park City Math. Ser.}, pages 3--159. Amer. Math. Soc.,
  Providence, RI, 1997.

\bibitem[Rog00]{rogers2000rank}
Nicholas~F Rogers.
\newblock Rank computations for the congruent number elliptic curves.
\newblock {\em Experimental Mathematics}, 9(4):591--594, 2000.

\bibitem[Rog04]{rogersthesis}
Nicholas~F. Rogers.
\newblock {\em Elliptic curves x(3) + y(3) = k with high rank}.
\newblock PhD thesis, 2004.
\newblock Accessed through ProQuest Dissertations and Theses; Last updated
  2020-05-06.

\bibitem[{Sag}20]{sage}
{Sage Developers}.
\newblock {\em {S}ageMath, the {S}age {M}athematics {S}oftware {S}ystem
  ({V}ersion 8.8)}, 2020.
\newblock {\tt https://www.sagemath.org}.

\bibitem[Tun83]{tunnell}
Jerrold~B. Tunnell.
\newblock A classical {D}iophantine problem and modular forms of weight
  {$3/2$}.
\newblock {\em Invent. Math.}, 72(2):323--334, 1983.
\newblock https://dx.doi.org//10.1007/BF01389327.

\bibitem[WDE{\etalchar{+}}14]{watkins2014ranks}
Mark Watkins, Stephen Donnelly, Noam~D Elkies, Tom Fisher, Andrew Granville,
  and Nicholas~F Rogers.
\newblock Ranks of quadratic twists of elliptic curves.
\newblock {\em Publications Math{\'e}matiques de Besan{\c{c}}on}, (2):63--98,
  2014.

\bibitem[Zag89]{zagier1989elliptische}
Don Zagier.
\newblock {\em Elliptische Kurven: Fortschritte und Anwendungen}.
\newblock Max-Planck-Institut f{\"u}r Mathematik, 1989.

\end{thebibliography}

\end{document}